\documentclass[11pt]{article}
\usepackage{amssymb}
\usepackage{color}
\usepackage{appendix}
\usepackage[normalem]{ulem}

\newcommand{\bN}{\mathbf{N}}

\newlength{\szer}

\newtheorem{defi}{Definition}
\newtheorem{nota}[defi]{Remark}

\newtheorem{teorema}[defi]{Theorem}

\newtheorem{propiedad}[defi]{Property}

\newenvironment{proof}[1][Proof]{\textbf{#1.} }{\
\rule{0.5em}{0.5em}}

\begin{document}
\title{On the Abhyankar-Moh inequality
\footnotetext{
     \noindent   \begin{minipage}[t]{4in}
       {\small
       2000 {\it Mathematics Subject Classification:\/} Primary 20M14;
       Secondary 32S05.\\
       Key words and phrases: semigroups, sequences of divisors, Abhyankar-Moh inequality.\\
       The second-named author was partially supported by the Spanish Projects PNMTM 2007-64007 and
       MTM2012-36917-C03-01.}
       \end{minipage}}}

\author{Roland D. Barrolleta, Evelia R.\ Garc\'{\i}a Barroso,  and Arkadiusz P\l oski}

\maketitle

\begin{abstract}
\noindent Abhyankar and  Moh in their fundamental paper on the embeddings of the line in the plane proved an important inequality which can be stated in terms of the semigroup associated with the branch at infinity of a plane algebraic curve. In this note we study the semigroups of integers satisfying the Abhyankar-Moh inequality and give 
a simple proof of the Abhyankar-Moh embedding theorem.
\end{abstract}

\noindent This short note is conceived as an appendix to \cite{GBP}. Our aim is to study the semigroup of integers appearing in connection with the Abhyankar-Moh inequality which is the main tool in proving the famous embedding line theorem (\cite{Abhyankar-Moh2}, Main Theorem). Since the Abhyankar-Moh inequality can be stated in terms of the semigroup associated with the branch at infinity of a plane algebraic curve (\cite{GBP}, Theorem 6.4) it is natural to consider the semigroups for which such an inequality holds. In what follows we use freely the properties of semigroups explained in \cite{GBP}, Preliminaries, pp. 3-6.

\medskip

\noindent Let $n>0$ be an integer. A semigroup $G\subseteq \bN$ will be called an {\em Abhyankar-Moh semigroup of degree} $n$ if the $n$-minimal sequence of generators $\overline{b}_0=n, \overline{b}_1,\dots, \overline{b}_h$ of $G$ satisfies the following conditions
\begin{enumerate}
\item [(G1)]Set $e_k=\gcd (\overline{b}_0,\ldots,\overline{b}_k)$  for $0\leq k \leq  h$ and $n_k=\frac{e_{k-1}}{e_k}$ for $1\leq k\leq h$. Then $e_h=1$ and $n_k>1$ for $1\leq k \leq  h$.

\item [(G2)] $n_{k-1}\overline{b}_{k-1}<\overline{b}_k$ for  $2\leq k\leq h$.

\item [(G3)] $n_{h-1}\overline{b_h}<n^2$.
\end{enumerate}

\noindent Our aim is to prove the following

\begin{teorema}
\label{th-AM} 
Let $G$ be an Abhyankar-Moh semigroup of degree $n>1$ and let $c$ be the conductor of $G$, that is the smallest element of $G$ such that all integers bigger than or equal to it are in $G$. Then $c\leq (n-1)(n-2)$. If $\overline{b}_0=n,\overline{b}_1,\ldots,\overline{b}_h$ is the $n$-minimal sequence of generators of $G$ then $c=(n-1)(n-2)$ if and only if $\overline{b}_k=\frac{n^2}{e_{k-1}}-e_k$
 for $k\in \{1,\ldots,h\}$.
\end{teorema}

\noindent \begin{proof}
Let $\delta_0=n$ and $\delta_k=\displaystyle\frac{n^2}{e_{k-1}}-\overline{b}_k $ for  $k\in \{1, \dots, h\}$. We have
$\delta_k=\displaystyle\frac{n^2-e_{k-1}\overline{b}_k}{e_{k-1}}\geq \displaystyle\frac{n^2-e_{h-1}\overline{b}_h}{e_{k-1}}>0$
and $\gcd(\delta_0,\ldots,\delta_k)=e_k$ for $k\in \{1, \dots, h\}$. Let $\gamma=\sum_{k=1}^h (n_k-1)\delta_k-\delta_0+1$. Since $\delta_k \geq e_k$ we get $\gamma \geq \sum_{k=1}^h (n_k-1)e_k-e_0+1=0$ with equality if and only if $\delta_k=e_k$ for $k\in \{0, \dots, h\}$.

\medskip

\noindent On the other hand 
\begin{eqnarray*}
\gamma&=&\sum_{k=1}^h (n_k-1)\left(\frac{n^2}{e_{k-1}}-\overline{b}_k\right)-n+1=\sum_{k=1}^h (n_k-1)\frac{n^2}{e_{k-1}}-n+1-\sum_{k=1}^h (n_k-1)\overline{b}_k\\
&=&(n-1)^2-\sum_{k=1}^h (n_k-1)\overline{b}_k=(n-1)(n-2)-c
\end{eqnarray*}

\noindent since $c=\sum_{k=1}^h (n_k-1)\overline{b}_k-\overline{b}_0+1$ by \cite{GBP}, Proposition 1.5, {\em 4}.

\medskip

\noindent Therefore $c\leq (n-1)(n-2)$ since $\gamma\geq 0$ and $c=(n-1)(n-2)$ if and only if $\gamma=0$ which is equivalent to $\delta_k=e_k$  that is $\displaystyle\frac{n^2}{e_{k-1}}-\overline{b}_k=e_k. $ 
\end{proof}

\medskip

\noindent Let  $n>1$ be an integer. A sequence of integers $e_0, \ldots, e_h$ will be called a {\em sequence of divisors of } $n$ if $e_k$ divides $e_{k-1}$ for $k\in \{1,\ldots,h\}$ and 
 $n=e_0>e_1>\cdots>e_{h-1}>e_h=1$.

\medskip

\noindent Now, we can give a simple description of the Abhyankar-Moh semigroups of degree $n>1$ with $c=(n-1)(n-2)$. For any sequence of divisors $e_0,e_1,\ldots,e_h$ of $n$, we put

\[G(e_0,\ldots,e_h)= \bN n +\bN (n-e_1) +\cdots+\bN \left(\frac{n^2}{e_{k-1}}-e_k\right)+\cdots+\bN \left(\frac{n^2}{e_{h-1}}-1\right).\]

\noindent Using Theorem \ref{th-AM} we check easily

\begin{propiedad}
\label{AM01}
A semigroup $G$ is an Abhyankar-Moh semigroup of degree $n>1$ with $c=(n-1)(n-2)$ if and only if $G=G(e_0,\ldots,e_h)$ where $n=e_0, e_1, \ldots, e_{h}$ is a sequence of divisors of $n$. 
\end{propiedad}

\noindent Here is another application of Theorem \ref{th-AM} 
\begin{propiedad}
\label{AM11}
Let $G$ be an Abhyankar-Moh semigroup of degree $n>1$ with $c=(n-1)(n-2)$ and let $n'=\min (G\backslash\{0\})$. Then $\gcd(n,n')=n-n'$. In particular, $n-n'$ divides $n$.
\end{propiedad}

\noindent \begin{proof}
Let $\overline{b}_0=n, \overline{b}_1,\ldots,\overline{b}_h$ be the $n$-minimal sequence of generators of $G$. The sequence $e_{k-1}\overline{b}_k$, 
$k\in \{1, \dots, h\}$ being increasing we get that $e_0\overline{b}_1\leq e_{h-1}\overline{b}_h<n^2$, that is $\overline{b}_1<n$. Consequently $\overline{b}_1$ is the smallest term of the sequence $\overline{b}_0,\overline{b}_1,\ldots,\overline{b}_h$ and $n'=\min (G\backslash\{0\})=\overline{b}_1$. By the second part of Theorem \ref{th-AM} $\overline{b}_1=e_1(n_1-1)=n-e_1$ and the property follows.
\end{proof}

\medskip

\begin{propiedad}
\label{AM21}
Let $G$ be an Abhyankar-Moh semigroup of degree $n$ with $c=(n-1)(n-2)$ and let $\overline{\beta_0},
\overline{\beta_1},\ldots, \overline{\beta_g}$ be the minimal system of generators of the semigroup $G$. Then
$n=\overline{\beta_1}$ or $n=2\overline{\beta_0}$.\\
\noindent If $n=\overline{\beta_1}$ then $G=\bN n+\bN (n-\epsilon_1) +\bN \left(\frac{n^2}{\epsilon_{1}}-\epsilon_1\right)+\cdots+\bN \left(\frac{n^2}{\epsilon_{g-1}}-1\right)$.\\
\noindent If $n=2\overline{\beta_0}$ then $G=\bN (n-\epsilon_0)+\bN \left(\frac{n^2}{\epsilon_0}-\epsilon_1\right) +\cdots+ \bN \left(\frac{n^2}{\epsilon_{g-1}}-1\right),$\\

\noindent where $\epsilon_k=\gcd\left(\overline{\beta_0},\overline{\beta_1},\ldots, \overline{\beta_k}\right)$ for $k\in \{0,\ldots,g\}$.
\end{propiedad}

\noindent \begin{proof}
The conditions (G1) and (G2) imply that $n\leq \overline{\beta_1}$ and $n\equiv 0$ (mod $\overline {\beta_0}$) if $n<{\beta_1}$. We claim that if $n\neq {\beta_1}$ then $n=2{\overline{\beta_0}}$. Indeed, if $n\neq {\beta_1}$ then $n=a{\overline{\beta_0}}$ for an integer $a>0$ and $n=b(n-\overline{\beta_0})$ for an integer $b>0$ by Property 
\ref{AM11}. Thus we get $a=(a-1)b$ which implies $a=2$. The second part of the property follows from \cite{GBP}, Proposition 1.5, {\em 3}.
\end{proof}

\begin{nota}
Suppose that the sequence $\overline{b_0}=n,\ldots,\overline{b_1},\ldots,\overline{b_h}$ satisfies the conditions \hbox{\rm (G1), (G2)} and \hbox{\rm (G3)}. Then the sequence of positive integers $\delta_0=\overline{b}_0=n$, $\delta_k=\displaystyle\frac{n^2}{e_{k-1}}-\overline{b}_k $ for  $k\in \{1, \dots, h\}$,  has the following properties

\begin{enumerate}
\item[(i)] The sequence of divisors $\gcd(\delta_0,\ldots,\delta_k)=e_k$, $k\in \{1, \dots, h\}$, is strictly decreasing, $e_h=1$.
\item[(ii)] $\delta_1<\delta_0$ and $\delta_k<n_{k-1}\delta_{k-1}$ for $k\in \{2, \dots, h\}$, where $n_k=\frac{e_{k-1}}{e_k}$.
\end{enumerate}

\noindent In contrast to the property $n_k\overline{b}_k \in \bN n+\cdots+\bN \overline{b}_{k-1}$ of the sequences $\overline{b}_0,\ldots,\overline{b}_h$ satisfying  the conditions \hbox{\rm (G1)} and \hbox{\rm (G2)}, in general, the condition $n_k\delta_k \in \bN n+\cdots+\bN \delta_{k-1}$ fails.
For example $(\overline{b}_0,\overline{b}_1,\overline{b}_2)=(6,2,17)$ is an Abhyankar-Moh sequence with the corresponding sequence $(\delta_0,\delta_1,\delta_2)=(6,4,1)$ but $n_2\delta_2=2\not\in \bN \delta_0+\bN \delta_1=\bN 6+\bN 4$.
\end{nota}

\noindent {\bf Abhyankar-Moh Embedding Line Theorem} (\cite{GBP}, Theorem 6.6).
\noindent Assume that $C$ is a rational projective irreducible curve of degree $n>1$ with one branch at infinity 
and such that the center of the branch at infinity $O$ is the unique singular point of $C$. Suppose that $C$ 
is permissible and let $n'$ be the multiplicity of $C$ at $O$. Then $n-n'$ divides $n$.

\medskip

\noindent \begin{proof}
It follows from \cite{GBP}, Theorem 6.4
that $\Gamma_O$ is an Abhyankar-Moh semigroup of degree $n$. 
Let $c$ be the conductor of the semigroup $\Gamma_O$. Using the Noether formula for the genus of projective plane curve we get $c=(n-1)(n-2)$. Then the theorem follows from Property \ref{AM11}.
\end{proof}

\medskip
\noindent {\small Roland David Barrolleta\\
Departament d'Enginyeria de la Informaci\'o i de les Comunicacions\\
Universitat Aut\`onoma de Barcelona\\
08193 - Bellaterra - Cerdanyola, Barcelona, Spain\\
e-mail:  rolanddavid.barrolleta@deic.uab.cat}

\medskip

\noindent {\small Evelia Rosa Garc\'{\i}a Barroso\\
Departamento de Matem\'aticas, Estad\'{\i}stica e I.O.\\
Secci\'on de Matem\'aticas, Universidad de La Laguna\\
38271 La Laguna, Tenerife, Espa\~na\\
e-mail: ergarcia@ull.es}

\medskip

\noindent {\small Arkadiusz P\l oski\\
Department of Mathematics\\
Kielce University of Technology \\
Al. 1000 L PP7\\
25-314 Kielce, Poland\\
e-mail: matap@tu.kielce.pl}

\end{document}